# HUYGENS WAVE EQUATIONS IN THE FIELD OF 2D-CWT

*Victor Vermehren Valenzuela*[1,2], *Hélio Magalhães de Oliveira*[2]

[1]Univ. Estadual do Amazonas, Manaus-AM, Brazil, [2]Univ.de Federal de Pernambuco, Recife-PE, Brazil

## ABSTRACT

**In this paper it is shown the performing of an optical transform to state the scalar diffraction in the formulation of the wavelet transform and the wave equations. From there, a bridge is build between equations of spherical waves presented in 1678 by Huygens and the continuous wavelet transform. For such a purpose, wavelets are introduced that meet the principles of waves and the properties of wavelets. The following equations are applied in solution to show a correspondence between the Huygens-Fresnel diffraction and the wavelet transform.**

*Index Terms*— Diffraction, 2D-CWT, Fresnel, wavelets, chirp.

## 1. INTRODUCTION

Wavelet and its Continuous and Discrete transformed (DWT and CWT, respectively) have emerged as a definitive tool in signal processing [1] and [2]. This is mainly due to the fact it provides an accurate local and global information identification of signals [3]. Another application experiencing great attention is the use of wavelet analysis techniques in optical processing [4-5], both treating images and diffraction, which extend wavelet analysis for two (2D) and three (3D) dimensions. However, an understanding of the relationship between the very essence of optics and the concept of wavelet, has been attempted in recent decades, through the implementation of optical devices or by representations in own domains [6-10]. Nevertheless, it had not yet been established with property the wave propagation carries inherent characteristics of the two-dimensional wavelet transform.

The scalar theory of diffraction involves the translation of the wave equation, which is a partial differential equation, in an integral equation. This can be used to analyze most imaging systems and phenomena within its area of validity [11]. The diffraction pattern of a 2-D object, $U(x, y)$, at a distance $z$ according to the theory of Fresnel diffraction region (approximate) can be written as

$$U_z(x, y) = \frac{e^{j\frac{2\pi z}{\lambda}}}{j\lambda z} \int\!\!\int_{-\infty}^{\infty} U(x_0, y_0)\, e^{j\frac{\pi}{\lambda z}[(x-x_0)^2 + (y-y_0)^2]}\, dx_0 dy_0 \quad (1)$$

The Equation (1) can be rewritten as

$$U_z(x, y) = U(x_0, y_0) ** h_z(x, y), \quad (2)$$

where $**$ denotes a 2-D convolution, $\lambda$ is the wavelength and the impulse response is expressed by

$$h_z(x, y) = \frac{e^{j\frac{2\pi z}{\lambda}}}{j\lambda z}\, e^{j\frac{\pi}{\lambda z}(x^2+y^2)}. \quad (3)$$

Equation (2) indicates that the Fresnel diffraction, which is based on the theory of waves of Huygens, exhibits characteristics of modern mathematical wavelet signals, namely a set of different diffraction patterns produces a set of "wavelet transformed images", from the same 2-D object.

The optical waves defined in Equation (3) are in fact the monochromatic spherical waves of Huygens. Under the Fresnel approximation, the Huygens-Fresnel equation can be interpreted as a wavelet transform. According Onural [6] it is straightforward to see from Equations (2) and (3) that the core of the diffraction convolution equation has indeed the properties shifting and scaling associated with a family of wavelet itself. This is shown by defining such a wavelet as follows:

$$h_z(x, y) = e^{j(x^2+y^2)}, \quad (4)$$

thus deriving from (5) the following family of wavelets

$$h_{abc}(x, y) = K_a\, \psi\!\left(\frac{x-b}{a}, \frac{y-c}{a}\right), \quad (5)$$

where $a = \sqrt{\lambda z / \pi}$, $b = \xi$, $c = \eta$, and the constant $K_a$ is

$$K_a = \frac{e^{j\frac{2\pi z}{\lambda}}}{j\lambda z} = e^{j\left\{\frac{2}{\pi a^2}\left[\left(\frac{\pi a}{\lambda}\right)^2 - \frac{\pi}{4}\right]\right\}}. \quad (6)$$

Note that the propagation distance $z$ can be interpreted as the wavelet scale factor.

Everything would be solved in the field of transform if the optical wavelet, as defined in the Equation (4), satisfies the conditions of wavelets, including the admissibility. Although the initial field equations $E(x_0, y_0, 0)$ establishes recovery via the diffracted field $E(x, y, z)$ by the inverse optical propagation, integrals do not converge to a solution even in FFT and IFFT computations, since they produce distortions in the framework of circular shapes [12].

Another challenge is that these wavelets are not located, neither in the space-time domain, nor in the frequency domain. Later, Onural himself called them chirp scaling functions instead of wavelets in the field of Fractional Fourier Transform [13].

## 2. THE PROPOSED WAVELET

It can be observed that the convolution in Equation (2) must be performed through a family of wavelet function for a continuous 2-D space that satisfies all conditions of waves. Classically, the wavelet functions are obtained from a function, indexed by two coefficients [1]:

$$\psi_{a,s}(x) = \frac{1}{\sqrt{s}} \psi\left(\frac{x-a}{s}\right), \quad (7)$$

where $a$ is a shift parameter, $s$ is a scale parameter, and $s^{1/2}$ is a normalization factor. The spread to higher dimensions consequently may be obtained and shown as follow,

$$\psi_{abs}(x, y) = \frac{1}{s} \psi\left(\frac{x-a}{s}, \frac{y-b}{s}\right). \quad (8)$$

From the point of view of wave equations and CWT, there had to be a family of wavelets that perform the same location or diffraction properties at the CWT. How to make a chirp function that tends to infinity become "behaved" (i.e., meet the requirements or properties of wavelets?)

The solution partly comes from a recent family of wavelets, the chirplets [14]. In part because although they have become bandwith limited functions due to Gaussian window applied, these wavelets are not even functions, since they have a term shift in exponential and increases with frequency and one end of chirp rate. They wavelets are building to specific image compression applications.

Then, taking into account the Equations (3-8) introduced in this paper, a new family of 2D wavelet called chirplet optics are presented below

$$\psi(x,y) = \frac{1}{\sqrt{\sqrt{2\pi}\sigma}} e^{-\frac{(x^2+y^2)}{4\sigma^2}} e^{j\frac{\pi}{\lambda z}(x^2+y^2)}, \quad (9)$$

where as usual $z$ is the propagation factor and $\lambda$ is the wavelength. Note that the value of $\sigma^2$ should be large enough, in order to take into account a maximum number of side lobes involved in the kernel. Furthermore, this wavelet becomes an even function as its analogous at (3) all-encompassing with the same representation of families of wavelets (5) and the continuous phase (6), limited bandwidth and holding the properties listed below.

The function $\psi(x,y)$ meets the requirements of wavelets of unity energy,

$$\int_{-\infty}^{\infty}\int_{-\infty}^{\infty} |\psi(x,y)|^2 \, dxdy = 1, \quad (10)$$

the requirement of zero mean

$$\int_{-\infty}^{\infty}\int_{-\infty}^{\infty} \psi(x,y) dxdy = 0, \quad (11)$$

and the nice property of vanishing moments

$$\int_{-\infty}^{\infty}\int_{-\infty}^{\infty} \psi(x,y) x^n y^n dxdy = 0, \quad (12)$$

with $n = 1, 2, 3...$ (computational verification up to 35).

The Fourier transform (FT) of the optical *chirplet* is expressed as follow:

$$\Psi(u,v) = \frac{\pi}{\alpha - j\beta} e^{-\frac{(u-w_0\cos\theta_u)^2 + (v-w_0\sin\theta_v)^2}{4(\alpha-j\beta)}}, \quad (13)$$

where $\alpha = 1/4\sigma^2$, $\beta = \pi/\lambda z$ and $w_0 = 2\pi/\lambda$. The optical chirplet function is defined as a cosenoidal complex chirp, whose amplitude is modulated by a Gaussian function and their frequency has quadratic sweep. Based on this model, the calculation of expression in the frequency domain occurs naturally through Fourier transform pair of a complex exponential signal quadratic, which dual is a complex Gaussian pulse too [15].

The graphs of chirplet optical wavelet and its corresponding FT are shown in Figure 1. Note the pulse characteristic of $\Psi(w)$ and the decreasing to zero value when it tends to the coordinates (0,0) and $\infty$, as expected.

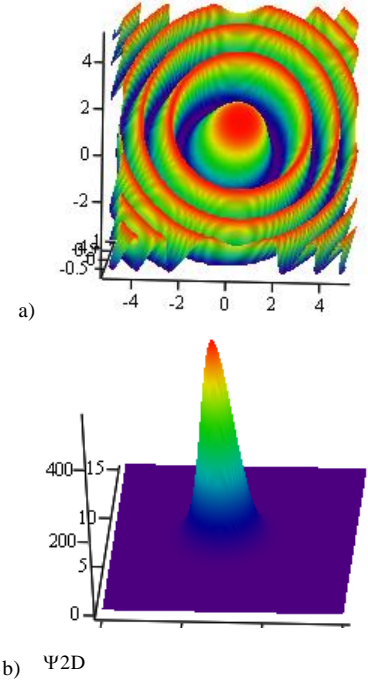

a)

b) $\Psi$2D

Figure 1. Graphs of the chirplet optical wavelet. a) 3-D wavelet representation. b) 3-D representation of Fourier transforma of the chirplet optical wavelet.

## 3. THE 2D-CWT

The continuous wavelet transform of a 1D signal $f(x)$ is defined as [16],

$$CWT(a,s) = \frac{1}{\sqrt{s}} \int_{-\infty}^{\infty} f(x) \psi_{a,s}^*(x) dx. \quad (14)$$

Note that Equation (14) takes the form of correlation between the input signal $f(x)$ and dilated and shifted mother wavelet $\psi_{as}(x)$, which is rewritten in the following convolution form [16]:

$$CWT(a,s) = (f * \psi_s^*)(a) = \int_{-\infty}^{\infty} \psi_s^*(a-x) f(x) dx. \quad (15)$$

Implementing the 2D-CWT, a pattern $g(x,y)$ is projected onto the wavelet $\psi_{a,b,s}$, by the translation of $a$ and $b$, in the $x$ and $y$ axes respectively, and scaling by $s$, the mother wavelet $\psi(x, y)$, as shown below.

$$CWT(a,b,s) = \frac{1}{s} \int_{-\infty}^{\infty}\int_{-\infty}^{\infty} g(x,y) \psi_{abs}^*(x,y) dxdy \quad (16)$$

or $\quad CWT(a,b,s) = (g ** \psi_s^*)(a,b).$

The wavelet transform, as in Equation (15) is therefore the convolution of the signal with a wavelet function. Thus, one can employ the convolution theorem to express the wavelet transform in terms of products of the Fourier transform of signal, $X(\omega)$, and the wavelet, $\Psi_{a,b}(\omega)$, as shown in [17]

$$CWT(a,s) = \int_{-\infty}^{\infty} X(\omega)\Psi_{a,s}^*(\omega)\,d\omega, \quad (17)$$

where the conjugate of the wavelet function is used, and $\omega=2\pi f$, whose the Fourier transform of the dilated and translated wavelet is

$$\Psi_{a,s}(\omega) = \frac{1}{\sqrt{s}} \int_{-\infty}^{\infty} \psi\left(\frac{t-a}{s}\right) e^{-j\omega t}\,dt. \quad (18)$$

Making the variable substitution $t' = (t-a)/s$, and so $dt = s\,dt'$, it is derived

$$\Psi_{a,s}(\omega) = \frac{1}{\sqrt{s}} \int_{-\infty}^{\infty} \psi(t') e^{-j\omega(st'+a)} s\,dt'. \quad (19)$$

Separating out the constant part of the exponential function and omitting the prime from the variable $t'$, one gets

$$\Psi_{a,s}(\omega) = \sqrt{s}\, e^{-j\omega a} \int_{-\infty}^{\infty} \psi(t) e^{-j(s\omega)t}\,dt. \quad (20)$$

The integral expression in the previous equation is just the FT of the wavelet at rescaled frequency $a\omega$. Consequently one can rewrite Equation (19) as:

$$\Psi_{a,s}(\omega) = \sqrt{s}\,\Psi(s\omega) e^{-j\omega a}. \quad (21)$$

The FT of the wavelet function conjugate is merely

$$\Psi_{a,s}^*(\omega) = \sqrt{s}\,\Psi^*(s\omega) e^{j\omega a}. \quad (22)$$

Equation (17) can thus be rewritten in expanded bidimensional form, yielding:

$$CWT_{FT}(a,b,s) = s \int\!\!\int_{-\infty}^{\infty} X(u,v)\Psi^*(au,bv) e^{j(au+bv)}\,du\,dv \quad (23)$$

## 4. THE 2D-ICWT

The diffraction involves reverse recovery image of an object to which the diffraction pattern was measured, for example, on a plane. In the case of the Fresnel and the Fraunhofer approximation inverting equations are obtained directly, only by inverting complexes operators [11]. In the case of very close fields, i.e., the one that approximates the illuminated object, it can be used angular spectrum representation techniques, but with the fulfillment of certain conditions. In both cases, the mathematical complexity is rather enormous, since it involves double integrals, solution typically are not easily convergent, as can be seen in the sequel.

a. Fraunhofer

$$U(x,y,z) = \frac{je^{-j\frac{2\pi}{\lambda}z_0}}{\lambda z_0} \int_{-\infty}^{\infty}\!\int_{-\infty}^{\infty} U(x_0,y_0,z_0) e^{-j\frac{\pi}{\lambda z_0}(x_0^2+y_0^2)}$$
$$\times e^{j\frac{2\pi}{\lambda z_0}(x_0 x + y_0 y)}\,dx_0\,dy_0 \quad (24)$$

b. Fresnel

$$U(x,y,z) = \frac{je^{-j\frac{2\pi}{\lambda}z_0}}{\lambda z_0} e^{-j\frac{\pi}{\lambda z_0}(x^2+y^2)} \int_{-\infty}^{\infty}\!\int_{-\infty}^{\infty} U(x_0,y_0,z_0)$$
$$\times e^{-j\frac{\pi}{\lambda z_0}(x_0^2+y_0^2)} e^{j\frac{2\pi}{\lambda z_0}(x_0 x + y_0 y)}\,dx_0\,dy_0 \quad (25)$$

c. Angular spectrum

$$U(x,y,z) = \int_{-\infty}^{\infty}\!\left[\int_{-\infty}^{\infty} U(x_0,y_0,z_0) e^{-j2\pi(f_x x_0 + f_y y_0)}\,dx_0\,dy_0\right]$$
$$\times e^{-jz_0\sqrt{\frac{4\pi^2}{\lambda^2}-4\pi^2(fx^2+fy^2)}} e^{j2\pi(f_x x + f_y y)}\,df_x\,df_y \quad (26)$$

where $z$ is the distance between the observation plain and the generation plain.

In contrast to the previous equations, there exists now the 2D-ICWT to reconstruct the input pattern through their wavelet chirplet optics decompositions:

$$ICWT(x,y) = \frac{1}{s^2 C} \int_{-\infty}^{\infty}\!\int_{-\infty}^{\infty} CWT(a,b,s)$$
$$\times K_s \psi\left(\frac{x-a}{s}, \frac{y-b}{s}\right) da\,db \quad (27)$$

where $s = \sqrt{\lambda z / \pi}$ denotes the observation scale parameter and $C$ is constant equal to (some sort of 2D admissibility condition).

$$C(u,v) = \int_{-\infty}^{\infty}\!\int_{-\infty}^{\infty} \frac{|\Psi(u,v)|^2}{|u\,v|}\,du\,dv < \infty. \quad (28)$$

Whereas $C$ must be a finite number so as to obtain the inverse wavelet transform, the wavelet chirplet optical meets this condition, as seen in Figure 1, and

$$\Psi(0,0) = 0. \quad (29)$$

The wavelet transform, as in Equation (27), is the convolution of the CWT with the wavelet function. Thus, one can also adopt the same procedure of Section 3 to express the 2D CWT products in terms of the Fourier transform of the CWT, $X(\omega)$, and wavelet, $\Psi_{a,b}(\omega)$, as shown [12]:

$$ICWT(x,y) = \frac{1}{sC} \int\!\!\int_{-\infty}^{\infty} CWT_{FT}(u,v)\Psi(su,sv) e^{-j(au+bv)}\,du\,dv \quad (30)$$

This is a result particularly useful in the same way as (23), since it assists the calculation of the inverse continuous wavelet transform into higher complex mathematical functions.

## 5. SIMULATION AND EXPERIMENTS

In this section some MathCad[TM] simulation are presented to corroborate the validity of the expressions of chirplet optical wavelet in the 2D-CWT, both under forward and reverse forms. Figure 2 shows the simulation the approximated Fresnel diffraction equation and optical chirplet 2D-CWT applied to a test object with a rectangular aperture 6 mm × 2 with unit amplitude, being illuminated by a uniform white light source.

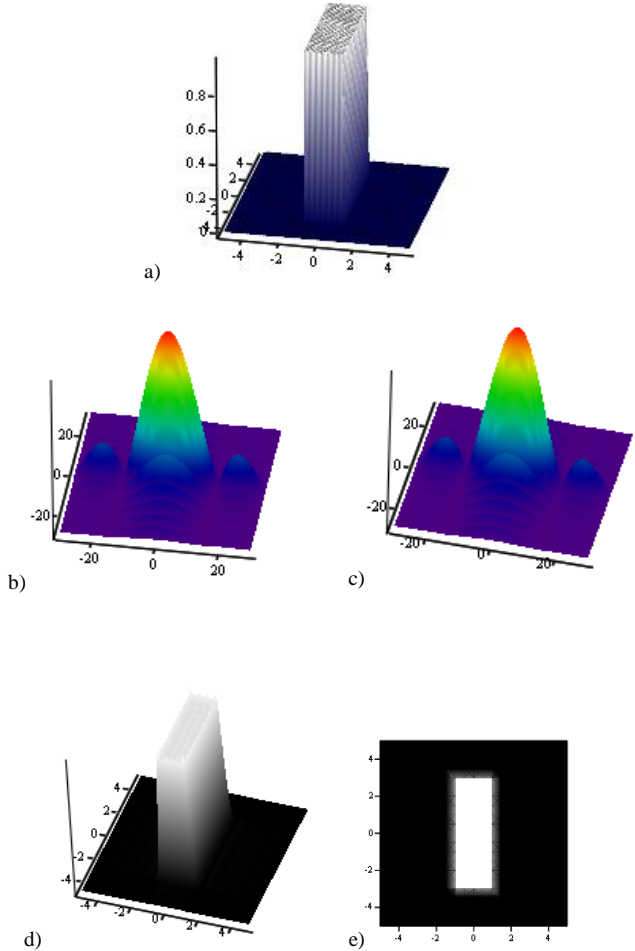

Figure 2. a) Object of the test at the generator plan. (b) Illustrates the application of the Fresnel diffraction equation at a distance z=54000mm and λ=550nm (c) Illustrates the amplitude pattern application of the chirplet 2D-CWT direct, at the same conditions as (b). (d) Shows the result of applying chirplet optical 2D-CWT reverse in the reconstruction of the object pattern generator and, (d) the mask that generates the amplitude pattern.

It is noteworthy that the diffraction reverse only the integrals of 2D-CWT converged to a solution (even the angular spectrum equations do not achieved any valuable result).

## 6. ON POTENTIAL APPLICATIONS OF THE NEW APPROACH

Despite the fact of knowing the foundations of the diffraction theory, well-established by Huygens-Fresnel for over 300 years, only in recent decades a renewed impulse appears with the introduction of its application in the fields of holography [19], tomography [20], ultrasound [21], and even in image compression [22].

These applications, especially in the biomedical area, present serious technological challenges and high-cost implications in the physical implementation of optical devices, as an attempt to overcome various barriers that cause errors and interferences.

Under this new formulation of the wavelet transform, with their characteristics of multiresolution and location in time or space-frequency, and while the inverse transform is a filter (13), now, e.g., can be applied to the 3D image retrieval in digital holographic microscopy [22]. In this application while eliminating the effect of the zero-order images, it is also possible to mitigate the blurry and twin image without spatial filters. The image of order-zero (DC) and the twin images are a pair of symmetrical terms in the zero-order spectrum. After that it is possible to adjust the filter in order to position the bandpass in range the required image (+1 order) as shown in Figure 3.

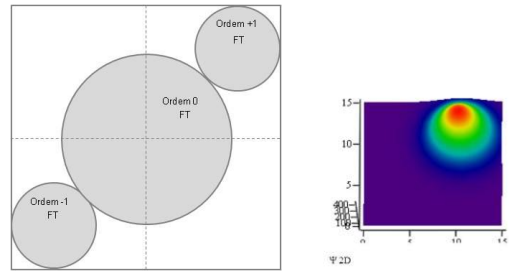

Figura 3. Filtragem com ICWT Chirplet Óptica. a) Representação das imagens no espectro. b) A transformada de Fourier da Wavelet Chirplet Óptica ajustada para filtragem da imagem ordem+1.

## 7. CONCLUDING REMARKS

It was successfully performed the wavelet optical chirplet, unprecedented for the 2D-analysis of waves, which builds a bridge between the equations of spherical waves introduced in 1678 by Christian Huygens and the continuous wavelet transform. Wide-ranging wavelets introduced, which meet the principles of waves and the properties of wavelets with the notorious fact of its reverse integral converge to a solution and its filter can be positioned in required range of application.

This development may product in significant reduction of costs, because only mathematical operations are performed in their implementation. An overview of preliminary experimental results illustrates this fact. At the same time the effect of zero-order term are eliminated, as well as the term obfuscation edge image and twin image without spatial filters holographic microscopy. A major advantage of the introduced technique is that it can be applied in the analysis of biological specimen's live dynamic images from digital holographic microscopy [23] more precisely. This analysis also can particularly be applied to microbial structures such as cancer cell, stem cell and botanical specimens.

The analysis of holograms to extract information related to objects can find a useful theoretical framework as a result of the relationship diffraction-wavelet presented. In addition, the noise filtering in holographic images is an extra and vast field of analysis.